\documentclass[smallextended,envcountsect]{svjour3}
\usepackage{graphicx,color,cite}
\usepackage{mathrsfs}
\usepackage{amsmath, amscd, amsfonts, amssymb, graphicx, color}
\usepackage[bookmarksnumbered, colorlinks, plainpages]{hyperref}
\usepackage{mathptmx}   
\usepackage[colorinlistoftodos]{todonotes}
\usepackage[left=4cm, right=4cm, top=3.5cm, bottom=3.5cm]{geometry}

\usepackage{latexsym}

\if{
\newtheorem{theorem}{Theorem}[section]

\newtheorem{lemma}{Lemma}[section]

\newtheorem{remark}{Remark}[section]

\newtheorem{corollary}{Corollary}[section]

}\fi

\renewcommand\theenumi{\roman{enumi}}

\newcounter{mycount}

\smartqed

\usepackage{etex}

\usepackage{mathtools}

\usepackage{enumitem}
\renewcommand\theenumi{(\roman{enumi})}
\renewcommand\theenumii{(\alph{enumii})}

\usepackage{csquotes}

\usepackage{todonotes}

\makeatletter
\let\orgdescriptionlabel\descriptionlabel
\renewcommand*{\descriptionlabel}[1]{
 \let\orglabel\label
 \let\label\@gobble
 \phantomsection
 \edef\@currentlabel{#1}
 \let\label\orglabel
 \orgdescriptionlabel{#1}
}
\def\th@plain{
 \thm@notefont{}
 \itshape
}
\def\th@definition{
 \thm@notefont{}
 \normalfont
}

\g@addto@macro\th@definition{\thm@headpunct{}}
\g@addto@macro\th@plain{\thm@headpunct{}}
\makeatother

\usepackage{subfig}

\usepackage[final]{showkeys}

\usepackage{etoolbox}

\usepackage{mathrsfs}

\usepackage[titletoc]{appendix}
\usepackage[doc]{optional}

\usepackage{soul}

\usepackage{colortbl,booktabs,
multirow}
\usepackage{xcolor}

\usepackage{cancel}

\usepackage{empheq}

\definecolor{myblue}{rgb}{.8, .8, 1}
\newcommand*\mybluebox[1]{
\colorbox{myblue}{\hspace{1em}#1\hspace{1em}}}

\usepackage{hyperref}
\hypersetup{
colorlinks=true,
linkcolor=blue,
citecolor=blue,
filecolor=magenta,
urlcolor=cyan
}

\usepackage[
open,
openlevel=2,
atend,
numbered
]{bookmark}

\usepackage[capitalize, nameinlink, noabbrev]{cleveref}
\crefname{equation}{}{}
\crefname{chapter}{Chapter}{Chapters}
\crefname{item}{item}{items}
\crefname{figure}{Figure}{Figures}
\crefname{theorem}{Theorem}{Theorems}
\crefname{lemma}{Lemma}{Lemmas}
\crefname{proposition}{Proposition}{Propositions}
\crefname{corollary}{Corollary}{Corollarys}
\crefname{definition}{Definition}{Definitions}
\crefname{fact}{Fact}{Facts}
\crefname{example}{Example}{Examples}
\crefname{algorithm}{Algorithm}{Algorithms}
\crefname{remark}{Remark}{Remarks}
\crefname{note}{Note}{Notes}
\crefname{notation}{Notation}{Notations}
\crefname{case}{Case}{Cases}
\crefname{exercise}{Exercise}{Exercises}
\crefname{question}{Question}{Questions}
\crefname{claim}{Claim}{Claims}
\crefname{enumi}{}{}

\usepackage{pgf}

\usepackage{array}
\usepackage{tabu}

\allowdisplaybreaks

\numberwithin{equation}{section}

\renewcommand\theenumi{(\roman{enumi})}
\renewcommand\theenumii{(\alph{enumii})}

\spnewtheorem*{Proof}{Proof.}{\bf}{\rm}
\begin{document}

\title{On Error Bounds of Inequalities in Asplund Spaces\thanks{Research  of the first author was supported by the National Natural Science Foundation of  China (grants 11971422 and 12171419), and funded by the Science and Technology Project of Hebei Education Department (No. ZD2022037) and the Natural Science Foundation of Hebei Province (A2022201002). Research of the second author  was partially supported  by  FASIC and a public grant as part of the Investissement d'avenir project, reference ANR-11-LABX-0056-LMH, LabEx LMH.}}

\titlerunning{On Error Bounds of Inequalities in Asplund Spaces}

\author{Zhou Wei  \and Michel Th\'era \and Jen-Chih Yao}

\institute{Zhou Wei\at Hebei Key Laboratory of Machine Learning and Computational Intelligence \& College of Mathematics and Information Science, Hebei University, Baoding, 071002, China\\ \email{weizhou@hbu.edu.cn}\\
Michel Th\'era \at XLIM UMR-CNRS 7252, Universit\'e de Limoges, Limoges, France\\\email{michel.thera@unilim.fr}\\
Jen-Chih Yao \at Research Center for Interneural Computing, China Medical University Hospital,
China Medical University, Taichung, Taiwan.\\ \email{yaojc@mail.cmu.edu.tw}
}

\date{Received: date / Accepted: date}
\dedication{Hedy Attouch passed away in October 2023 and the aim of this contribution is 
to celebrate  his memory as mathematician, as a colleague and as a friend.  Hedy's  personality,   enthusiasm,   cordiality,   passion for research will remain an example for the community.}
\maketitle

\begin{abstract}
Error bounds are central objects in optimization theory and its applications. They were for a long time restricted only to the theory before becoming  over the course of time a field of itself.  This paper is devoted to the study  of error bounds of a general inequality defined by a proper lower semicontinuous function on an Asplund space. Even though the results of the dual characterization on the error bounds of a  general inequality (if one drops the convexity assumption) may not be valid, several necessary dual conditions are still obtained in terms of  Fr\'echet/Mordukhovich subdifferentials   of the concerned function at points in the solution set.
%
Moreover, for an inequality  defined  by a composite-convex function  that is    to say by a function which is  the composition of a convex function with a smooth mapping, such dual conditions also turn out  to be sufficient to have the error bound property. Our work is an extension of  results on dual characterizations of convex inequalities to the possibly non-convex case.



\end{abstract}

\section{Introduction}

The main goal  in this paper is to study error bounds of a general inequality in Asplund spaces. The starting point of theory of error bounds can be traced back to the pioneering work by Hoffman \cite{28} for systems of affine functions in which this author proved that for a given matrix and a vector, the distance from the candidate point to the solution set is bounded  from above by some scalar constant times the norm of the residual error. Hoffman's result was extensively and intensively studied by many authors. Robinson \cite{54} extended Hoffman's work to a system of convex inequalities which defines a bounded feasible region with a nonempty interior. Subsequently, Mangasarian \cite{42} derived a global error bound for
differentiable convex inequalities under the Slater condition and an  asymptotic constraint qualification. Then, Auslender and Crouzeix \cite{1} used closed proper convex functions satisfying an interior condition on their domains to improve Mangassarian's work for systems which are not necessary differentiable on the whole space.  Pang \cite{52} gave a comprehensive survey of an extensive theory and rich applications of error bounds for inequalities and optimization systems and solution sets to equilibrium problems. Lewis and Pang \cite{41} gave a systematic and unified treatment for the existence of a global error bound for a convex inequality system. For the summary of the theory of error bounds and their various applications, the readers are invited to consult bibliographies \cite{AT2022,3,7,18,24,29,ioffe-book,34,jourani,Kruger-LY,Luke,44,penot-book,55} and references therein for more details.

Error bounds have been well-recognized in optimization and variational analysis and played an important role in various aspects like the sensitivity analysis of linear programs (cf. \cite{Rob73,Rob77}), the  convergence analysis of descent methods for linearly constrained minimization (cf. \cite{Gul92,IuD90,TsL92,TsB93}), the convex feasibility problem of finding a point in the intersection of a finite collection of closed convex sets (cf. \cite{5,6,7}) and  the domain of image reconstruction (cf. \cite{16}). Also,  error bounds  have been proved to closely relate with weak sharp minima of functions,  metric regularity/subregularity, as well as Aubin property/calmness of set-valued mappings (cf. \cite{3,BD1,BD2,14,ioffe-JAMS-1,ioffe-JAMS-2,54,59,60} and references therein).

To study error bounds, a large literature is devoted to provide dual characterizations and criteria in terms of subdifferentials. It is known that two directions    have been followed to describe error bounds. One direction makes use of the subdifferentials of the given function at points outside the solution set. To the best of our knowledge, \cite{Ioffe1979} was among  of the first papers of such kind to state sufficient conditions for error bounds of a constraint system in terms of  the Clarke subdifferential. In 2009,  Ngai and Th\'era \cite{NT2009} discussed error bounds for systems of lower semicontinuous functions in Asplund spaces and provided sufficient conditions in terms of Fr\'echet subdifferentials of the concerned functions at points outside the solution set. The other direction  focuses on   subdifferentials of the given function at points inside the solution set. Along this direction, Lewis and Pang \cite{41} studied error bounds for convex inequality systems and provided necessary conditions via subdifferentials and normal cones.  In 2003,  Ngai and Th\'era \cite{NT2004} provided an error bound estimate and an implicit multifunction theorem in terms of smooth subdifferentials and abstract subdifferentials. In 2004,  Zheng and Ng \cite{ZN2004} proved dual characterizations of error bounds for convex inequalities in terms of subdifferentials and normal cones. In 2010, Fabian, Henrion, Kruger and Outrata \cite{FHKO2010} presented a general classification scheme of necessary and sufficient criteria for error bounds and used several derivative-like objects both from the primal as well as from the dual space to characterize error bounds of extended-real-valued functions on a Banach space.

Inspired by \cite{FHKO2010,NT2009,ZN2004}, we  aim to study error bounds of a general inequality defined by a proper semicontinuous function on an Asplund space and  in addition   to give dual criteria ensuring error bounds. Although the dual characterization results on error bounds for general inequality may not hold (if dropping the convexity assumption), several necessary conditions for error bounds are still obtained  in terms of Fr\'echet/Mordukhovich subdifferentials of the given function at points in the solution set. Further, such conditions are proved to be also sufficient for error bounds of the composite-convex inequality (defined by a composition of a convex function with a smooth mapping). Our work is an extension of dual characterizations for error bounds via subdifferentials to the possible non-convex case. 

The   paper is organized as follows. Several preliminaries  and known results will be given in Section 2. Section 3 is devoted to error bounds of a general inequality in an Asplund space. Several dual necessary conditions for error bounds are obtained in terms of Fr\'echet/limiting subdifferentials of the concerned function at points in the solution set (see \cref{th3.1,th3.2,th3.3}). Then we further discuss error bounds of a composite-convex inequality (defined by a composition of a convex function with a smooth mapping) and prove that  such dual conditions are also sufficient for error bounds in this case (see \cref{th3.4,th3.5}). Conclusions of this paper are presented in Section 4.

 \medskip

\setcounter{equation}{0}

\section{Preliminaries}
Let $\mathbb{X}$ be a Banach space  whose norm is denoted by $\Vert\cdot\Vert$.  $\mathbb{X}^*$ stands for the (continuous) dual of $\mathbb{X}$ equipped with the dual norm $\Vert\cdot\Vert_*,$ while    $ \mathbb{B}_{\mathbb{X}}$ 
and  $ \mathbb{B}_{\mathbb{X}^*}$ are respectively   the closed unit ball  in $\mathbb{X}$  and $\mathbb{X}^*$ and $\mathbb{B}(x,r)= \{y\in \mathbb{X} : \Vert x-y\Vert <r\}$  is the open ball      ball centered at $x$  with radius equal to   $r>0$.
 Given  a multifunction $\Phi: \mathbb{X}\rightrightarrows \mathbb{X}^*$, the symbol
\begin{equation*}
\begin{array}r
\mathop{\rm Limsup}\limits_{y\rightarrow x}\Phi(x):=\Big\{x^*\in \mathbb{X}^*: \exists \ {\rm sequences} \ x_n\rightarrow x \ {\rm and} \ x_n^*\stackrel{w^*}\longrightarrow x^* \ {\rm with}\  \\
x_n^*\in \Phi(x_n) \ {\rm for\ all \ } n\in \mathbb{N} \Big\}
\end{array}
\end{equation*}
signifies the  sequential Painlev\'{e}-Kuratowski outer/upper limit of $\Phi(x)$ as $y\rightarrow x$, where the notation $x_n^*\stackrel{w^*}\longrightarrow x^*$ means that the sequence $\{x_n^*\}$ converges to $x^*$ with respect to the weak$^*$-topology.

For a subset $\Omega$ of $\mathbb{X}$, we denote by ${\rm cl}(\Omega)$ and ${\rm bd}(\Omega)$ the closure and the boundary of $\Omega$, respectively. For any $x\in \mathbb{X}$, the distance from $x$ to  the subset $\Omega$ is given by $$\mathbf{d}(x, \Omega):=\inf\{\|x-u\|:u\in\Omega\}.$$ 

Let $A$ be a closed set and $a\in A$. For $\varepsilon\geq 0$, we denote
$$
\hat N_{\varepsilon}(A, x):=\left\{x^*\in
\mathbb{X}^*:\limsup\limits_{x\stackrel{A}\rightarrow a}\frac{\langle x^*,
x-a\rangle}{\|x-a\|}\leq \varepsilon\right\}
$$
the set of $\varepsilon$-\textit{normals} to $A$ at $a$, where $x\stackrel{A}\rightarrow a$ means that $x$ tends to $a$ while staying in $A$. When $\varepsilon=0$,  we will write  $\hat N(A,a)$ instead of $\hat N_{0}(A,a)$. $\hat N(A,x)$  is a closed convex cone  called the {\it Fr\'echet normal cone} to $A$ at $a$.

The  {\it Mordukhovich (or limiting) normal cone of $A$ at $a$  is the set  $N(A, a)$ which is defined by}
$$
N(A, a):=\mathop{\rm Limsup}_{x\stackrel A\rightarrow a, \varepsilon\downarrow
0}\hat N_{\varepsilon}(A, x).
$$
It is known that
$$\hat N(A,a)\subset N(A, a).$$
If $A$ is convex, the Fr\'echet normal cone coincides with the Mordukhovich normal cone and reduces to the normal cone in the sense of convex analysis; that is,
$$
\hat{N}(A,a)=N(A, a)=\{x^*\in \mathbb{X}^*:\;\langle
x^*,x-a\rangle\leq 0\;\; \forall x\in A\}.
$$

\begin{empheq}[box =\mybluebox]{equation}\label{2.1}
N(A, a)=\mathop{\rm Limsup}_{x\stackrel A\rightarrow
a}\hat N(A, x)
\end{empheq}
Thus $x^*\in N(A, a)$ if and only if there exist $x_n\stackrel{A}\rightarrow a$ and $x^*_n\stackrel{w^*}\rightarrow x^*$ such that $x_n^*\in \hat N(A, x_n)$ for all $n$.

Let  $f :\mathbb{X}\rightarrow \mathbb{R}\cup\{+\infty\}$ be a proper
lower semicontinuous function and $x\in
\mathrm{dom\,}f:=\{y\in \mathbb{X}: f(y)<+\infty\}$. We denote by
$${\rm epi\,}f:=\{(x, r)\in \mathbb{X}\times \mathbb{R}: f(x)\leq r\}
$$ 
the {\it epigraph} of $f$. We denote by 
\begin{eqnarray*}
\hat{\partial}f(x):=\{x^*\in \mathbb{X}^* : (x^*, -1)\in \hat N({\rm epi\,}f, (x, f(x)))\}
\end{eqnarray*}
the {\it Fr\'{e}chet subdifferential} of $f$ at $x$ and   by
\begin{eqnarray*}
\partial f(x):=\{x^*\in \mathbb{X}^* : (x^*, -1)\in N({\rm epi\,}f, (x, f(x)))\},
\end{eqnarray*}
the {\it Mordukhovich subdifferential} of $f$ at $x$, respectively.


We denote by
$$
\hat\partial^{\infty}f(x):=\{x^*\in \mathbb{X}^*: (x^*,0)\in \hat N({\rm epi\,}f,(x,f(x)))\}
$$
the {\it Fr\'echet singular subdifferential} of $f$ at $x$ and   by
$$
\partial^{\infty}f(x):=\{x^*\in \mathbb{X}^*: (x^*,0)\in N({\rm epi\,}f,(x,f(x)))\}
$$
the {\it Mordukhovich singular subdifferential} of $f$ at $x$., respectively. When $f$ is Lipschitzian around $x$, the Fr\'echet and Mordukhovich singular subdifferentials of $f$ at $x$ are  trivial; that is, $\hat\partial^{\infty}f(x)=\partial^{\infty}f(x)=\{0\}$.

Well known are the facts that  $$
\hat{\partial}f(x)\subseteq\partial f(x),
$$
and
$$
\hat{\partial}f(x)=\left\{x^*\in \mathbb{X}^*:\;\liminf\limits_{u\rightarrow x}
\frac{f(u)-f(x)-\langle x^*,u-x\rangle}{\|u-x\|}\geq0\right\}.
$$

Recall that a Banach space $\mathbb{X}$ is said to be an Asplund space if every continuous convex real-valued  function defined on a nonempty open convex subset $D\subseteq \mathbb{X}$ is Fr\'echet differentiable at each point of some dense $G_{\delta}$ subset of $D$.  There is a plethora of properties equivalent to Asplundness and the  readers are invited to consult \cite{Ph} or \cite{Yost93} for the definition of Asplund spaces and their equivalences. 
It is well known that among all these  equivalences,  
 $\mathbb{X}$ is Asplund    if and only if every separable subspace  of $\mathbb{X}$ has a separable dual. In particular, every reflexive Banach  space is an Asplund space.


Variational analysis in  the framework of Asplund  spaces    has been deeply carried out  by  Mordukhovich and  Shao \cite{  MS1996} and developed in several books  and papers, see \cite{Mordukhovich}  and the recent book by Thibault \cite [Comment page 555] {thibault}. Two  of the most important results in \cite{  MS1996}  are these following formulae that will be a key ingredient in our study: 
\vskip 2mm
The first one:
\begin{empheq}[box =\mybluebox]{equation}\label{2.1}
N(A, a)=\mathop{\rm Limsup}_{x\stackrel A\rightarrow
a}\hat N(A, x)
\end{empheq}
Equivalently,  $x^*\in N(A, a)$ if and only if there exist $x_n\stackrel{A}\rightarrow a$ and $x^*_n\stackrel{w^*}\rightarrow x^*$ such that $x_n^*\in \hat N(A, x_n)$ for all $n$.
\vskip 2mm
The second one:
\begin{empheq}[box =\mybluebox]{equation}\label{2.3}
\partial f(x)=\mathop{\rm Limsup}_{y\stackrel f\rightarrow
x}\hat{\partial}f(y)\  {\rm and} \  \partial^{\infty} f(x)=\mathop{\rm Limsup}_{y\stackrel f\rightarrow
x,\lambda\downarrow 0}\lambda\hat{\partial}f(y),
\end{empheq}
where  the notation $y\stackrel f\rightarrow x$ means that $y\rightarrow x$ with $f(y)\rightarrow f(x)$. \vskip 2mm
\noindent Equivalently,   $x^*\in\partial f(x)$ if and only if there exist $x_n\stackrel{f}\rightarrow x$ and $x_n^*\stackrel{w^*}\rightarrow x^*$ such that $x_n^*\in\hat{\partial}f(x_n)$ for all $n$, and $x^*\in\partial^{\infty} f(x)$ if and only if there exist $x_n\stackrel{f}\rightarrow x$, $\lambda_n\downarrow 0$ and $x_n^*\in\hat{\partial}f(x_n)$ such that $\lambda_nx_n^*\stackrel{w^*}\rightarrow x^*$.

\medskip

The following important fuzzy calculus result was proved by Fabian \cite{Fabian1989}. Readers could also consult \cite[Theorem 3]{Fabian1989} and \cite[Theorem 2.33]{Mordukhovich} for the details.

\begin{lemma}\label{lem2.1}
Let $\mathbb{X}$ be an Asplund space and $f_i:\mathbb{X}\rightarrow \mathbb{R}\cup\{+\infty\}$ $(i=1,2)$ be proper lower semicontinuous functions such that one of them  is Lipschitz continuous around $x\in{\it dom\,}f_1\cap {\it dom}\,f_2$. Then for any $x^*\in\hat\partial (f_1+f_2)(x)$ and $\varepsilon>0$, there exist $x_i\in \mathbb{B}(x,\varepsilon)$ with $|f_i(x_i)-f_i(x)|<\varepsilon \;(i=1,2)$ such that
$$
x^*\in \hat\partial f_1(x_1)+\hat\partial f_2(x_2)+\varepsilon \mathbb{B}_{\mathbb{X}^*}.
$$
\end{lemma}

The following Lemma 2.2 (see \cite[Lemma 3.6]{NT2001}  and also  \cite[Lemma 3.7]{AusDanThi05}) and Lemma 2.3 (see \cite[Theorem 3.1]{ZN2008}) will be used in our analysis later.
\begin{lemma}\label{lem2.2}
Let $\mathbb{X}$ be an Asplund space and $A$ be a nonempty closed subset of
$\, \mathbb{X}$. Let $x\in \mathbb{X}$ and $x^*\in \hat{\partial}\mathbf{d}(\cdot,
A)(x)$. Then for any $\varepsilon>0$ there exist $a\in A$ and
$a^*\in \hat N(A, a)$ such that
$$
\|x-a\|<\mathbf{d}(x,
A)+\varepsilon \ \ {\it and} \ \ \|x^*-a^*\|<\varepsilon.
$$
\end{lemma}


\begin{lemma}\label{lem2.3}
Let $\mathbb{X}$ be  an Asplund space,  $A$ be a nonempty closed
subset of $\,\mathbb{X}$ and $x\not\in A$. Then for any $\beta\in (0,\;1)$
there exist $z\in {\rm bd}(A)$ and $z^\ast \in
\hat{N}(A, z)$ with $\|z^\ast \|=1$ such that
$$\beta\|x-z\|<\min\{\mathbf{d}(x, A),\;\langle z^\ast ,x-z\rangle\}.$$
\end{lemma}

The following lemma is cited from \cite[Theorem 3.5.8]{BorweinZhu}.

\begin{lemma}\label{lem2.4}
Let $\mathbb{X}$ be Fr\'echet smooth Banach space and $f_i:\mathbb{X}\rightarrow \mathbb{R}\cup\{+\infty\}$ $(i=1,\cdots,m)$ be proper lower semicontinuous functions and $f:=\max\{f_1,\cdots, f_m\}$. Suppose that $x^*\in\hat \partial f(\bar x)$. Then for any $\epsilon>0$ and any weak$^*$ neighborhood $U^*$ of $\,0$ in $\mathbb{X}^*$, there exist $x_k\in \mathbb{B}(\bar x,\epsilon)$ with $|f_k(\bar x)-f_k(x_k)|<\epsilon$, $x_k^*\in\hat\partial f_k(x_k)$ and $\lambda_k\geq 0(k=1,\cdots,m)$ such that
\begin{equation}\label{2.3}
  \left|\sum_{k=1}^{m}\lambda_k-1\right|<\epsilon\ \ {\it and} \ \ x^*\in\sum_{k=1}^{m}\lambda_kx_k^*+U^*.
\end{equation}
\end{lemma}

We conclude this section with the following lemma which is of independent interest.

\begin{lemma}\label{lem2.5}
Let $f:\mathbb{X}\rightarrow \mathbb{R}\cup\{+\infty\}$ be a proper lower semicontinuous function and $(z,r)\in{\rm epi\,}f$. Then
\begin{equation}\label{2.4}
  \hat N({\rm epi\,}f, (z,r))\subseteq\hat N({\rm epi\,}f, (z,f(z))).
\end{equation}
Assume further that $f$ is continuous at $z$. Then
\begin{equation}\label{2.5}
N({\rm epi\,}f, (z,r))\subseteq N({\rm epi\,}f, (z,f(z))).
\end{equation}
\end{lemma}

{\bf Proof.} Let $(x^*,t)\in\hat N({\rm epi\,}f, (z,r))$. Then
$$
\limsup_{(x,\alpha)\stackrel{{\rm epi\,}f}\longrightarrow (z,r)}\frac{\langle(x^*,t), (x-z,\alpha-r)\rangle}{\|(x-z,\alpha-r)\|}\leq 0.
$$
Then for any $(y,\beta)\stackrel{{\rm epi\,}f}\longrightarrow (z,f(z))$, one has $(y,\beta+r-f(z))\stackrel{{\rm epi\,}f}\longrightarrow (z,r)$ as $f(z)\leq r$ and thus
$$
\limsup_{(y,\beta)\stackrel{{\rm epi\,}f}\longrightarrow (z,f(z))}\frac{\langle(x^*,t), (y-z,(\beta+r-f(z))-r)\rangle}{\|(x-z,\beta-f(z))\|}\leq 0.
$$
This implies that $(x^*,t)\in\hat N({\rm epi\,}f,  (z,f(z)))$ and thus \eqref{2.4} holds. 

Let $(x^*,t)\in  N({\rm epi\,}f,  (z,r))$. Then there exist $\epsilon_k\downarrow 0, (z_k,r_k)\stackrel{{\rm epi\,}f}\longrightarrow (z,r)$ and $(x_k^*,t_k)\stackrel{w^*}\longrightarrow (x^*,t)$ such that $(x_k^*,t_k)\in\hat N_{\epsilon_k}({\rm epi\,}f, (z_k,r_k))$ for all $k$. Then for any $k$, one has
\begin{equation}\label{2.6}
  \limsup_{(x,\alpha)\stackrel{{\rm epi\,}f}\longrightarrow (z_k,r_k)}\frac{\langle(x_k^*,t_k), (x-z_k,\alpha-r_k)\rangle}{\|(x-z_k,\alpha-r_k)\|}\leq \epsilon_k.
\end{equation}
Let $(y,\beta)\stackrel{{\rm epi\,}f}\longrightarrow (z_k,f(z_k))$. Then $(y,\beta+r_k-f(z_k))\stackrel{{\rm epi\,}f}\longrightarrow (z_k,r_k)$ as $f(z_k)\leq r_k$ and it follows from \eqref{2.6} that
$$
\limsup_{(y,\beta)\stackrel{{\rm epi\,}f}\longrightarrow (z_k,f(z_k))}\frac{\langle(x_k^*,t_k), (y-z,(\beta+r-f(z))-r)\rangle}{\|(x-z_k,\beta-f(z_k))\|}\leq \epsilon_k.
$$
This means that $(x_k^*,t_k)\in\hat N_{\epsilon_k}({\rm epi\,}f, (z_k, f(z_k)))$. Noting that $f$ is continuous at $z$, it follows that $(x^*,t)\in N({\rm epi\,}f, (z,f(z)))$ as $z_k\rightarrow z$. Hence \eqref{2.5} holds. The proof is complete. \hfill$\Box$

\begin{remark}
It should be noted that  \eqref{2.5} may not hold necessarily if dropping the continuity assumption. For example, let $\mathbb{X}:=\mathbb{R}$ and $f:\mathbb{X}\rightarrow\mathbb{R}$ be defined by
\begin{equation*}
  f(x):=\left\{
  \begin{aligned}
  1, \ & x>0,\\
  x, \ & x\leq 0.
  \end{aligned}
  \right.
\end{equation*}
Let $(z,r):=(0,1)$. Then one can verify that $f$ is not continuous at $z$ and 
$$
N({\rm epi\,}f,(z,f(z)))=\{(t,s):t+s\geq 0, s\leq 0\}.
$$
Further one has
$$
(0,-1)\in N({\rm epi\,}f,(z,r))\backslash N({\rm epi\,}f,(z,f(z))),
$$
which implies that \eqref{2.5} does not hold. 
\end{remark}


\setcounter{equation}{0}
\section{Main Results}
In this section, we study error bounds of a general inequality defined by a proper lower semicontinuous function on an Asplund space and  we aim to give dual criteria for ensuring error bounds. As main results, it is proved that several necessary dual conditions for error bounds are provided in terms of Fr\'echet/Mordukhovich subdifferentials of the concerned function at points in the solution set. Further, for a composite-convex inequality defined by a composition of a convex function with a smooth mapping, such dual conditions are proved to be also sufficient for the existence of error bounds. We begin with the definition of error bounds of a general inequality.

\medskip
{\it Throughout of this section, unless stated otherwise, we always suppose that $\mathbb{X}$ is an Asplund space and $f:\mathbb{X}\rightarrow \mathbb{R}\cup\{+\infty\}$ is a proper lower semicontinuous function.}\\

We consider the following inequality: 
\begin{empheq}[box =\mybluebox]{equation}\label{3.1}
  \text{Find }\; x\in \mathbb{X}\;  \text{such that}\;f(x)\leq 0.
\end{empheq}
We denote by $\mathbf{S}_f:=\{x\in \mathbb{X}: f(x)\leq 0\}$  the solution set of inequality \eqref{3.1}.


Recall that inequality \eqref{3.1} is said to have a {\it local error bound} at $\bar x\in \mathbf{S}_f$ if there exist $\tau,\delta\in (0, +\infty)$ such that
\begin{equation}\label{4.8a}
  \mathbf{d}(x,\mathbf{S}_f)\leq\tau f_+(x)\ \ \forall x\in \mathbb{B}(\bar x,\delta),
\end{equation}
where $f_+(x):=\max\{f(x), 0\}$. The scalar
\begin{equation}\label{3.2a}
  \tau(f,\bar x):=\inf\{\tau>0:\exists \ \delta>0 \; \text{s.t. }\ \eqref{4.8a}\ {\rm holds}\}
\end{equation}
is called the {\it modulus of local error bound} of $f$ at $\bar x$.

The following theorem gives a necessary condition for local error bounds of inequality \eqref{3.1} in terms of Mordukhovich subdifferentials and singular subdifferentials.
\begin{theorem}\label{th3.1}
  Suppose that inequality \eqref{3.1} has a local error bound at $\bar x\in \mathbf{S}_f$ and $f$ is continuous at $\bar x$. Then 
\begin{equation}\label{3.3}
  N(\mathbf{S}_f, \bar x)\subseteq [0,+\infty)\partial f(\bar x)+\partial^{\infty} f(\bar x).
\end{equation}
\end{theorem}

{\it Proof.} Since inequality \eqref{3.1} has a local error bound at $\bar x$, then there exist $\tau,\delta\in (0,+\infty)$ such that \eqref{4.8a} holds. For any $(x, r)\in \mathbb{X}\times
\mathbb{R}$, let $\|(x, r)\|_{\tau}:=\frac{1}{\tau}\|x\|+|r|$. Then $(\mathbb{X}\times \mathbb{R}, \|\cdot\|_{\tau})$ is a Banach space and the unit ball of the
dual space of $(\mathbb{X}\times \mathbb{R}, \|\cdot\|_{\tau})$ is
$(\frac{1}{\tau}\mathbb{B}_{\mathbb{X}^*})\times \mathbb{B}_{\mathbb{R}}$. Let $\delta_0:=\frac{\delta}{2}$. We claim that
\begin{equation}\label{4.11}
\mathbf{d}(x, \mathbf{S}_f)\leq \tau(\mathbf{d}_{\|\cdot\|_{\tau}}((x, r), {\rm epi\,}f)+|r|)\ \ \forall (x, r)\in \mathbb{B}(\bar x, \delta_0)\times \mathbb{R}.
\end{equation}

Indeed, suppose on the contrary that there exists $(x_0, r_0)\in \mathbb{B}(\bar x, \delta_0)\times \mathbb{R}$ such that
$$
\mathbf{d}(x_0, \mathbf{S}_f)> \tau(\mathbf{d}_{\|\cdot\|_{\tau}}((x_0, r_0),
{\rm epi\,}f)+|r_0|).
$$
Then there exists $(u,r)\in {\rm epi\,}f$ such that
$$
\mathbf{d}(x_0, \mathbf{S}_f)> \tau(\frac{1}{\tau}\|u-x_0\|+|r_0-r|+|r_0|),
$$
and thus
$$
\mathbf{d}(x_0, \mathbf{S}_f)>\|u-x_0\|+\tau |r|.
$$
Note that $f(u)\leq r$ and one can verify that
\begin{equation}\label{4.12}
  \mathbf{d}(x_0, \mathbf{S}_f)>\|u-x_0\|+\tau f_+(u).
\end{equation}
Since
$$
\|u-\bar x\|\leq\|u-x_0\|+\|x_0-\bar x\|<\mathbf{d}(x_0, \mathbf{S}_f)+\|x_0-\bar x\|\leq
2\|x_0-\bar x\|<\delta,
$$
it follows from \eqref{4.8a} and \eqref{4.12} that
$$
  \mathbf{d}(x_0, \mathbf{S}_f)>\|u-x_0\|+\tau f_+(u)\geq \|u-x_0\|+\mathbf{d}(u,\mathbf{S}_f)\geq \mathbf{d}(x_0, \mathbf{S}_f),
$$
and  a contradiction. Hence \eqref{4.11} holds.

Let $x^*\in N(\mathbf{S}_f, \bar x)$. Then there exist $x_n\stackrel{\mathbf{S}_f}\rightarrow \bar x$ and $x_n^*\stackrel{w^*}\rightarrow x^*$ such that $x_n^*\in \hat N(\mathbf{S}_f, x_n)$. Applying the  Banach-Steinhaus Theorem, we can assume that $\|x_n^*\|\leq M$ for some $M>0$ and all $n$. Then for any $n$ sufficiently large, by \cite[Corollary 1.96]{Mordukhovich}, one has
$$
\frac{x^*_n}{M}\in\hat N(\mathbf{S}_f,x_n)\cap \mathbb{B}_{\mathbb{X}^*}=\hat\partial \mathbf{d}(\cdot, \mathbf{S}_f)(x_n)
$$
and thus for any $\varepsilon>0$ there exists $\delta_1\in (0,\delta_0-\|x_n-\bar x\|)$ such that
\begin{equation}\label{4.13}
 \Big \langle \frac{x^*_n}{M}, u-x_n\Big\rangle\leq \mathbf{d}(u, \mathbf{S}_f)+\tau\varepsilon\|u-x_n\|\ \ \forall u\in \mathbb{B}(x_n,\delta_1).
\end{equation}
Since $\mathbb{B}(x_n,\delta_1)\subseteq \mathbb{B}(\bar x,\delta_0)$, it follows from \eqref{4.11} and \eqref{4.13} that
$$
   \Big \langle \frac{x^*_n}{M}, u-x_n\Big\rangle\leq \tau(\mathbf{d}_{\|\cdot\|_{\tau}}((u,r), {\rm epi\,}f)+|r|)+\tau\varepsilon\|u-x_n\|\ \ \forall (u,r)\in \mathbb{B}(x_n,\delta_1)\times \mathbb{R}.
$$
Let $\Phi(u,r):=\mathbf{d}_{\|\cdot\|_{\tau}}((u, r), {\rm epi\,}f)+|r|.$ Then one has
\begin{equation}\label{4.14}
  \Big(\frac{x_n^*}{\tau Mx}, 0 \Big)\in\hat\partial \Phi(x_n, 0).
\end{equation}
By \cref{lem2.1}, there exist $(z_n,r_n)\in \mathbb{B}(x_n,\frac{1}{n})\times(-\frac{1}{n}, \frac{1}{n})$ and  $s_n\in (-\frac{1}{n}, \frac{1}{n})$ such that
$$
(\frac{x_n^*}{\tau M}, 0)\in\hat\partial \mathbf{d}_{\|\cdot\|_{\tau}}((\cdot,\cdot), {\rm epi\,}f)(z_n,r_n)+\{0\}\times [-1,1]+\frac{1}{n}\mathbb{B}_{\mathbb{X}^*}\times [-\frac{1}{n}, \frac{1}{n}].
$$
Thus there exist $(z^*_n,-\alpha_n)\in \hat\partial \mathbf{d}_{\|\cdot\|_{\tau}}((\cdot,\cdot), {\rm epi\,}f)(z_n,r_n)$ and $\gamma_n\in [-1,1]$ such that
\begin{equation}\label{4.15}
\Big(\frac{x_n^*}{\tau M}, 0\Big)\in(z^*_n,-\alpha_n)+(0,\gamma_n)+\frac{1}{n}\mathbb{B}_{\mathbb{X}^*}\times [-\frac{1}{n}, \frac{1}{n}].
\end{equation}
By virtue of \cref{lem2.2}, there exist $(\tilde{z}_n,\tilde{r}_n)\in{\rm epi\,}f$ and $(\tilde{z}^*_n, -\tilde{\alpha}_n)\in\hat N({\rm epi\,}f, (\tilde{z}_n,\tilde{r}_n))$ such that $$ \|(\tilde{z}_n,\tilde{r}_n)-(z_n,r_n)\|<\mathbf{d}_{\|\cdot\|_{\tau}}((z_n,r_n), {\rm epi\,}f)+\frac{1}{n}$$ and
\begin{equation}\label{4.16}
 \|(\tilde{z}_n^*,-\tilde{\alpha}_n)-(z_n^*,-\alpha_n)\|<\frac{1}{n}.
\end{equation}
Then $\tilde{\alpha}_n\geq 0$ for all $n$. Since $\mathbb{X}$ is an Asplund space, we know  that $\mathbb{B}_{\mathbb{X}^*}$ is weak$^*$ sequentially compact. Note that
$$
(z^*_n,\alpha_n)\in \hat\partial \mathbf{d}_{\|\cdot\|_{\tau}}((\cdot,\cdot), {\rm epi\,}f)(z_n,r_n)\subseteq\frac{1}{\tau}\mathbb{B}_{\mathbb{X}^*}\times [-1,1],
$$
and thus, considering subsequences if necessary, we can assume that
$$z_n^*\stackrel{w^*}\rightarrow z^*\in\frac{1}{\tau}\mathbb{B}_{\mathbb{X}^*}, \ \alpha_n\rightarrow \alpha\in [0,1], \ r_n\rightarrow r \in [-1, 1] \ {\rm and} \ \gamma_n\rightarrow\gamma\in [-1, 1] . $$
 By \eqref{4.16}, one has
$$
\tilde{z}_n^*\stackrel{w^*}\rightarrow z^*, \ \tilde{r}_n\rightarrow r\ {\rm and} \ \tilde{\alpha}_n\rightarrow \alpha.
$$
Since $f$ is continuous at $\bar x$, one has
$$
f(\bar x)=\lim_{n\rightarrow \infty}f(\tilde{z}_n)\leq \lim_{n\rightarrow \infty}\tilde{r}_n=r
$$
and it follows from \cref{lem2.5} that
\begin{equation}\label{4.17}
  (z^*,-\alpha)\in N({\rm epi\,}f, (\bar x, r))\subseteq N({\rm epi\,}f, (\bar x, f(\bar x))).
\end{equation}
By taking limits with respect to  the weak$^*$-topology in \eqref{4.15} as $n\rightarrow \infty$, one has
\begin{equation}\label{4.21a}
 x^*=\tau M z^*.
\end{equation}
From \eqref{4.17}, one has $\alpha\geq 0$. If $\alpha=0$, then $z^*\in\partial^{\infty}f(\bar x)$ and so $x\in \partial^{\infty}f(\bar x)$, which implies that \eqref{3.3} holds. If $\alpha>0$, then $\frac{z^*}{\alpha}\in\partial f(\bar x)$ and
$$
x^*=\tau M z^*\in [0,\tau M\alpha]\partial f(\bar x)\subseteq [0, +\infty)\partial f(\bar x).
$$
This means that \eqref{3.3} holds. The proof is complete.\hfill$\Box$

\medskip

In terms of Fr\'echet subdifferentials and singular subdifferentials, the following theorem gives a sharper necessary condition for a local error bound  for nequality \eqref{3.1}.

\medskip
\begin{theorem}\label{th3.2}
Suppose that inequality \eqref{3.1} has a local error bound at $\bar x\in \mathbf{S}_f$ and set 
$$\mathbb{B}_f(x,\epsilon):=\{u\in \mathbb{B}(x,\epsilon):|f(u)-f(x)|<\epsilon\}.$$ 
Then there exist $\tau,\delta_0>0$ such that for any $\epsilon>0$, 
the inclusion \begin{equation}\label{3.12}
  \hat N(\mathbf{S}_f,x)\cap \mathbb{B}_{\mathbb{X}^*}\subseteq \bigcup\left\{[0,(1+\epsilon)\tau]\hat\partial f(u)+\hat\partial^{\infty} f(v):u,v\in \mathbb{B}_f(x,\epsilon)\right\}+\epsilon \mathbb{B}_{\mathbb{X}^*}
\end{equation}
holds for all $x\in \mathbf{S}_f\cap \mathbb{B}(\bar x,\delta_0)$. 
\end{theorem}

{\it Proof.} By \eqref{4.8a}, there exist $\tau,\delta\in (0,+\infty)$ such that \eqref{4.8a} holds. Using the proof of Theorem 3.1, by defining $\|(x, r)\|_{\tau}:=\frac{1}{\tau}\|x\|+|r|$ for any $(x,r)\in \mathbb{X}\times \mathbb{R}$, one has that \eqref{4.11} holds with $\tau>0$ and $\delta_0=\frac{\delta}{2}>0$. Let $x\in \mathbb{B}(\bar x, \delta_0)\cap \mathbf{S}_f$ and $x^*\in \hat N(\mathbf{S}_f, x)\cap \mathbb{B}_{\mathbb{X}^*}$. 
Then  $x^*\in \hat\partial \mathbf{d}(\cdot, \mathbf{S}_f)(x)$ 
 by \cite[Corollary 1.96]{Mordukhovich}. Let $\epsilon>0$. Take $\epsilon_1>0$ such that 
$$
\max\{2\tau\epsilon_1, 2\epsilon_1, \frac{1}{\tau}\epsilon_1+\epsilon_1\}<\epsilon.
$$
Let $\Phi(u,r):=\mathbf{d}_{\|\cdot\|_{\tau}}((u, r), {\rm epi\,}f)+|r|$. Similarly to the proof of Theorem 3.1, one can verify that
\begin{equation}\label{4.19}
 \Big (\frac{x^*}{\tau}, 0\Big)\in\hat\partial \Phi(x, 0).
\end{equation}
By \cref{lem2.3}, there exist $(z_1,r_1)\in \mathbb{B}(x,\epsilon_1)\times(-\epsilon_1, \epsilon_1)$ and $s_1\in (-\epsilon_1, \epsilon_1)$ such that
$$
\Big(\frac{x^*}{\tau}, 0\Big)\in\hat \partial \mathbf{d}_{\|\cdot\|_{\tau}}((\cdot,\cdot), {\rm epi\,}f)(z_1,r_1)+\{0\}\times [-1,1]+\epsilon_1\mathbb{B}_{\mathbb{X}^*}\times [-\epsilon_1, \epsilon_1].
$$
Thus there exist $(z^*_1,-\alpha_1)\in \hat\partial \mathbf{d}_{\|\cdot\|_{\tau}}((\cdot,\cdot), {\rm epi\,}f)(z_1,r_1)$ and $\gamma_1\in [-1,1]$ such that
\begin{equation}\label{4.20}
\Big(\frac{x^*}{\tau}, 0\Big)\in(z^*_1,-\alpha_1)+(0,\gamma_1)+\epsilon_1\mathbb{B}_{\mathbb{X}^*}\times [-\epsilon_1, \epsilon_1].
\end{equation}
By virtue of \cref{lem2.2}, there are $(\tilde{z}_1,\tilde{r}_1)\in{\rm epi\,}f$ and $(\tilde{z}^*_1, -\tilde{\alpha}_1)\in\hat N({\rm epi\,}f, (\tilde{z}_1,\tilde{r}_1))$ such that $\|(\tilde{z}_1,\tilde{r}_1)-(z_1,r_1)\|<\mathbf{d}_{\|\cdot\|_{\tau}}((z_1,r_1), {\rm epi\,}f)+\epsilon_1$ and
\begin{equation}\label{4.21}
\|(\tilde{z}_1^*,-\tilde{\alpha}_1)-(z_1^*,-\alpha_1)\|<\epsilon_1.
\end{equation}
Then $(\tilde{z}^*_1, -\tilde{\alpha}_1)\in\hat N({\rm epi\,}f, (\tilde{z}_1, f(\tilde{z}_1)))$ by \cref{lem2.5} and
\begin{equation*}
\|(\tilde{z}_1, \tilde{r}_1)- (x,0))\|\leq\|(\tilde{z}_1, \tilde{r}_1)-(z_1, r_1)\|+\|(z_1, r_1)- (x,0))\|<\frac{1}{\tau}\epsilon_1+\epsilon_1.
\end{equation*}
From $(\tilde{z}^*_1, -\tilde{\alpha}_1)\in\hat N({\rm epi\,}f, (\tilde{z}_1, f(\tilde{z}_1)))$, one has $\tilde{\alpha}_1\geq 0$. If $\tilde{\alpha}_1=0$, then $\tilde{z}^*_1\in\hat\partial^{\infty}f(\tilde{z}_1)$. By virtue of \eqref{4.20} and \eqref{4.21}, one has
$$
\|\frac{x^*}{\tau}-\tilde{z}^*_1\|\leq\|\frac{x^*}{\tau}-{z}^*_1\|+\|{z}^*_1-\tilde{z}^*_1\|<2\epsilon_1
$$
and consequently
$$
x^*\in\hat\partial^{\infty}f(\tilde{z}_1)+2\tau\epsilon_1 \mathbb{B}_{\mathbb{X}^*}\subseteq\hat\partial^{\infty}f(\tilde{z}_1)+ \epsilon \mathbb{B}_{\mathbb{X}^*}
$$
(thanks to the choice of $\epsilon_1$). This implies that \eqref{3.12} holds.

If $\tilde{\alpha}_1>0$, then $\frac{\tilde{z}^*_1}{\tilde{\alpha}_1}\in\hat\partial f(\tilde{z}_1)$ and $\tilde{z}^*_1\in[0,1+2\epsilon_1]\hat\partial f(\tilde{z}_1)$ by \eqref{4.21}. By virtue of \eqref{4.20} and \eqref{4.21}, one has
$$
x^*\in [0, (1+2\epsilon_1)\tau]\hat\partial f(\tilde{z}_1)+2\tau\epsilon_1 \mathbb{B}_{\mathbb{X}^*}\subseteq [0, (1+\epsilon)\tau]\hat\partial f(\tilde{z}_1)+\epsilon \mathbb{B}_{\mathbb{X}^*}.
$$
Hence \eqref{3.12} holds. The proof is complete. \hfill$\Box$\\

\medskip

In \cref{th3.2}, a necessary condition for the etrror bound of inequality \eqref{3.1} was established in terms of Fr\'echet subdifferentials of the given function and normal cones at those points inside the solution set in \eqref{3.12}. This gives rise to a natural question whether $\epsilon$ in \eqref{3.12} can be taken as $0$. We have no answer to this question. However, we are now in a position to consider the circumstance when $\epsilon\downarrow 0$. The following theorem shows that this can be done in   finite-dimensional spaces.

\begin{theorem}\label{th3.3}
Suppose that $\mathbb{X}$ is finite-dimensional, inequality \eqref{3.1} has a local error bound at $\bar x\in \mathbf{S}_f$ and that $f$ is continuous at $\bar x$. Then there exists $\tau>0$ such that 
\begin{equation}\label{3.16}
  N(\mathbf{S}_f, \bar x)\cap \mathbb{B}_{\mathbb{X}^*}\subseteq [0, \tau]\partial f(\bar x)+ \partial^{\infty} f(\bar x).
\end{equation}
\end{theorem}

{\bf Proof.}  According to the local error bound property at $\bar x$, there exist $\tau,\delta\in (0,+\infty)$ such that \eqref{4.8a} holds. Let $x^*\in N(\mathbf{S}_f, \bar x)\cap \mathbb{B}_{\mathbb{X}^*}$ with $\|x^*\|>0$. Since $\mathbb{X}$ is finite-dimensional, there exist $x_n\stackrel{\mathbf{S}_f}\rightarrow \bar x$ and $x_n^*\stackrel{\|\cdot\|}\rightarrow x^*$ such that $x_n^*\in \hat N(\mathbf{S}_f, x_n)$. Then for any $n$ sufficiently large, one has $\frac{x_n^*}{\|x_n^*\|}\in \hat \partial \mathbf{d}(\cdot, \mathbf{S}_f)(x_n)$ by \cite[Corollary 1.96]{Mordukhovich} and thus for any $\epsilon>0$ there is $\delta_n\in (0,\delta-\|\bar x-x_n\|)$ such that
$$
\Big\langle\frac{x_n^*}{\|x_n^*\|},u-x_n\Big\rangle\leq \mathbf{d}(u, \mathbf{S}_f)+\tau\epsilon\|u-x_n\|, \ \ \forall u\in \mathbb{B}(x_n,\delta_n).
$$
This and \eqref{4.8a} imply that for any $u\in \mathbb{B}(x_n,\delta_n)$, one has
$$
\Big\langle\frac{x_n^*}{\|x_n^*\|},u-x_n\Big\rangle\leq \tau f_+(u)+\tau\epsilon\|u-x_n\|
$$
and thus
$$
\frac{x_n^*}{\tau\|x_n^*\|}\in\hat\partial f_+(x_n).
$$
(thanks to $f_+(x_n)=0$). Since $\mathbb{X}$ is finite-dimensional, it follows from \cref{lem2.4} that there exist $\lambda_n\in [0,1+\frac{1}{n}]$, $u_n\in \mathbb{B}(x_n,\frac{1}{n})$ with $|f(u_n)-f(x_n)|<\frac{1}{n}$ and $u_n^*\in\hat\partial f(u_n)$ such that
\begin{equation}\label{3-17}
 \left\| \frac{x_n^*}{\tau\|x_n^*\|}-\lambda_nu_n^*\right\|<\frac{1}{n}.
\end{equation}
Without loss of generality, considering a  subsequence if necessary, we can assume that $\lambda_n\rightarrow \lambda\in [0, 1]$. We divide $\lambda$ into two cases:
\begin{itemize}
  \item [$(a)$] $\lambda=0$. Note that $f$ is continuous at $\bar x$ and then $u_n\rightarrow \bar x$ with $f(u_n)\rightarrow f(\bar x)$. By letting $n\rightarrow \infty$ in \eqref{3-17}, one has 
      $$
      \frac{x^*}{\tau\|x^*\|}\in\partial^{\infty}f(\bar x)
      $$ 
      and thus \eqref{3.16} holds.
  \item [$(b)$] $\lambda>0$.  By letting $n\rightarrow \infty$ in \eqref{3-17}, one has
  $$
  u_n^*\rightarrow\frac{x^*}{\lambda\tau\|x^*\|}\in\partial f(\bar x)
  $$
  (as $f$ is continuous at $\bar x$). This means that 
  $$
  x^*\in\lambda\tau\|x^*\|\partial f(\bar x)\subseteq [0,\tau]\partial f(\bar x)
  $$
  (thanks to $\lambda\leq 1$ and $\|x^*\|\leq 1$). Hence \eqref{3.16} holds.
\end{itemize}

 The proof is complete. \hfill$\Box$

\medskip

\begin{remark}\label{rem3.1}
 It is known from \cref{th3.3} that the validity of \eqref{3.16} is necessary for the local error bound. However, this condition may not be sufficient for the local error bound even in the finite-dimensional space.  For example, let $\mathbb{X}:=\mathbb{R}^2$, 
 $$
 A_1:=(\mathbb{R}_-\times\mathbb{R}_+)\cup\{(x_1, x_2)\in \mathbb{R}^2:x_1+x_2=0, x_1\geq 0\}
 $$ 
and
$$
A_2:=\left\{(x_1, x_2)\in\mathbb{R}^2_+: (x_1-1)^2+x_2^2\leq 1\ \ {\rm and} \ \ x_1^2+(x_2-1)^2\leq 1\right\}.
$$
Define $f(x):=\max\{\mathbf{d}(x,A_1),\mathbf{d}(x,A_2)\}$ for all $x\in X$, and take $\bar x=(0,0), e=(\frac{\sqrt{2}}{2},\frac{\sqrt{2}}{2})$. Then $\mathbf{S}_f=A_1\cap A_2=\{\bar x\}$ and by computation, one has
$$
\partial \mathbf{d}(\cdot,A_1)(\bar x)\supseteq ([0, 1]\times\{0\})\cup(\{0\}\times[-1,0])\cup\{e\}
 $$
and
$$
\partial \mathbf{d}(\cdot,A_2)(\bar x)=N(A_2,\bar x)\cap \mathbb{B}_{\mathbb{X}^*}=\left\{(x_1, x_2):x_1^2+x_2^2\leq 1,x_1\leq 0,x_2\leq 0\right\}.
$$
Then one can verify that there exists $\varepsilon>0$ such that
$$
\partial f(\bar x)\supseteq  \varepsilon \mathbb{B}_{\mathbb{X}^*}=N(\mathbf{S}_f,\bar x)\cap \varepsilon \mathbb{B}_{\mathbb{X}^*},
$$
which implies that  \eqref{3.16} holds for $x=\bar x$.

On the other hand, if $x^{(k)}=\big(\frac{1}{k}, (\frac{2}{k}-\frac{1}{k^2})^{\frac{1}{2}}\big)$ for each $k\in \mathbb{N}$, one has
$$
\mathbf{d}(x^{(k)}, \mathbf{S}_f)=\big(\frac{2}{k}\big)^{\frac{1}{2}}, \mathbf{d}(x^{(k)}, A_1)=\frac{1}{k}\,\,\,{\rm and}\,\,\,\mathbf{d}(x^{(k)}, A_2)=0,
$$
and thus
$$
\frac{\mathbf{d}(x^{(k)}, \mathbf{S}_f)}{f(x^{(k)})}= (2k)^{\frac{1}{2}}\rightarrow \infty.
$$
This means that inequality $f(x)\leq 0$ does not have a  local error bound at $\bar x$.
\end{remark}

\medskip

Next, we consider error bounds for a convex-composite inequality  defined by a composite-convex function (i.e. the  composition  of a convex function with a continuously differentiable mapping), and  we use Fr\'echet subdifferentials and normal cones to prove dual characterizations of error bounds.

\medskip

 {\it Throughout the remainder of this section, we suppose that $\mathbb{Y}$ is an Asplund space, $\psi:\mathbb{X}\rightarrow \mathbb{Y}$ is continuously differentiable and $g:\mathbb{Y}\rightarrow \mathbb{R}\cup\{+\infty\}$ is a proper semicontinuous convex function on $\mathbb{Y}$}.

\medskip

Let $f:=g\circ\psi$. We consider the following convex-composite   inequality:
\begin{empheq}[box =\mybluebox]{equation}\label{3.21}
  (g\circ\psi)(x)\leq 0.
\end{empheq}
We denote by $\mathcal{S}$ the solution set of \eqref{3.21}, and denote by
\begin{equation}\label{3.22}
  \mathbf{S}_{g}:=\{y\in \mathbb{Y}: g(y)\leq 0\}
\end{equation} 
the solution set of  the convex inequality $g(y)\leq 0$.

\medskip

Now, we are in a position to present dual characterizations for error bounds of composite-convex inequality \eqref{3.21} when $\psi$  is smooth.

\begin{theorem}\label{th3.4}
Let $f:=g\circ\psi$ be such that ${\rm bd}(\mathbf{S}_g)\subseteq g^{-1}(0)$ and $\bar x\in  \mathcal{S}$ such that $\triangledown \psi(\bar x)$ is surjective. Then convex-composite inequality \eqref{3.21} has a local error bound at $\bar x$ if and only if there exist $\tau,\delta\in (0,+\infty)$ such that
\begin{equation}\label{4.19}
  \hat N(\mathcal{S}, x)\cap \mathbb{B}_{\mathbb{X}^*}\subseteq [0,\tau]\hat\partial f(x)
\end{equation}
holds for all $x\in \mathcal{S}\cap \mathbb{B}(\bar x, \delta)$.
\end{theorem}

{\bf Proof.} Since ${\rm bd}(\mathbf{S}_g)\subseteq g^{-1}(0)$, then it follows that 
\begin{equation}\label{3.24}
  {\rm bd}(\mathcal{S})\subseteq \psi^{-1}({\rm bd}(\mathbf{S}_g))\subseteq \psi^{-1}(g^{-1}(0))=f^{-1}(0).
\end{equation}
Note that $\mathcal{S}=\psi^{-1}(\mathbf{S}_g)$ and by applying \cite[Lemma 4.2]{WTY2023}, there are $\kappa, L>0$ and $r_0>0$ such that 
\begin{equation}\label{4.30}
\hat N(\mathcal{S},x)\cap \kappa \mathbb{B}_{\mathbb{X}^*}\subseteq \triangledown \psi(x)^*(N(\mathbf{S}_g, \psi(x))\cap \mathbb{B}_{\mathbb{Y}^*})\subseteq \hat N(\mathcal{S}, x)\cap L\mathbb{B}_{\mathbb{X}^*}
\end{equation}
holds for all $x\in \mathbb{B}(\bar x,r_0)\cap \mathcal{S}$.

{\it The necessity part}. By the local error bound, there exist $\tau>0$ and $\delta \in (0,r_0)$ such that \eqref{4.8a} holds. Let $x\in \mathcal{S}\cap \mathbb{B}(\bar x, \delta)$ and $x^*\in \hat N(\mathcal{S},x)\cap \mathbb{B}_{\mathbb{X}^*}$. Then $x^*\in \hat\partial \mathbf{d}(\cdot, \mathcal{S})(x)$ by \cite[Corollary 1.96]{Mordukhovich} and thus for any $\epsilon>0$, there is $\delta_1\in (0, \delta-\|x-\bar x\|)$ such that
\begin{equation*}
  \langle x^*,u-x\rangle\leq \mathbf{d}(u,\mathcal{S})+\tau\epsilon\|u-x\|,\ \ \forall u\in \mathbb{B}(x,\delta_1).
\end{equation*}
By using \eqref{4.8a}, one has
\begin{equation*}
  \langle x^*, u-x\rangle\leq \tau f_+(u)+\tau\epsilon\|u-x\|,\ \ \forall u\in \mathbb{B}(x,\delta_1).
\end{equation*}
This implies that
$$
\frac{x^*}{\tau}\in \hat\partial f_+(x)
$$
and it follows from \cite[Theorem 3.36 and Theorem 3.46]{Mordukhovich} that
$$
\frac{x^*}{\tau}\in [0,1]\partial f(x),
$$
By virtue of \cite[Lemma 4.1]{WTY2023}, one has $\partial f(x)=\hat\partial f(x)$ and thus
$$
x^*\in [0,\tau]\hat\partial f(x),
$$
which implies that \eqref{4.19} holds.

{\it The sufficiency part}.  
Let $\epsilon>0$ be such that $\epsilon<\kappa$. Since $\psi$ is continuously differentiable at $\bar x$, there exists $r_1 >0$ such that $r_1<\min\{r_0, \delta\}$ and
\begin{equation}\label{4.21a}
\|\psi(x)-\psi(u)-\triangledown \psi(u)(x-u)\|<\epsilon\|x-u\|,\ \forall x,u\in \mathbb{B}(\bar x,r_1).
\end{equation}
Take $\delta_1:=\frac{r_1}{2}$ and let $x\in \mathbb{B}(\bar x,\delta_1)\backslash \mathcal{S}$. Then $\mathbf{d}(x, \mathcal{S})\leq\|x-\bar x\|<\delta_1$. Choose $\beta\in (0,1)$ such that
$\beta>\max\{\frac{\mathbf{d}(x, \mathcal{S})}{\delta_1}, \frac{\epsilon}{\kappa}\}$. By virtue of \cref{lem2.3}, there exist $a\in {\rm bd}(\mathcal{S})$ and $a^*\in\hat N(\mathcal{S},a)$ with $\|a^*\|=1$ such that
\begin{equation}\label{4.22}
  \beta\|x-a\|\leq\min\{\mathbf{d}(x, \mathcal{S}), \langle a^*, x-a\rangle\}.
\end{equation}
Then
$$
\|a-\bar x\|\leq\|a-x\|+\|x-\bar x\|<\frac{\mathbf{d}(x, \mathcal{S})}{\beta}+\delta_1<2\delta_1=r_1
$$
and it follows from \eqref{4.19} and \cite[Lemma 4.1]{WTY2023} that there exist $\lambda\in[0,\tau], y^*\in\partial g(\psi(a))$ such that
\begin{equation*}\label{4.23}
  a^*= \triangledown\psi(a)^*(\lambda y^*).
\end{equation*}
By \eqref{4.30}, there is $y_1^*\in N(\mathbf{S}_g,\psi(a))\cap \mathbb{B}_{\mathbb{Y}^*}$ such that 
$$
\kappa a^*=\triangledown\psi(a)^*(y_1^*).
$$
Since $\triangledown\psi(a)^*$ is one-to-one, then one has $\lambda y^*=\frac{1}{\kappa}y_1^*$ and it follows from \eqref{4.21a} that
\begin{eqnarray*}
&&\langle a^*, x-a\rangle= \lambda\langle\triangledown \psi(a)^*(y^*),x-a\rangle= \langle \lambda y^*,\triangledown \psi(a)(x-a)\rangle\\
&=&\langle \lambda y^*,\psi(x)-\psi(a)\rangle+\langle \lambda y^*,\triangledown \psi(a)(x-a)+\psi(a)-\psi(x)\rangle\\
&\leq& \lambda\big(g(\psi(x))-g(\psi(a))\big)+\frac{1}{\kappa}\|\psi(x)-\psi(a)-\triangledown \psi(a)(x-a))\|\\
&\leq& \tau(f(x)-f(a))+\frac{\epsilon}{l}\|x-a\|\\
&=&\tau f_+(x)+\frac{\epsilon}{\kappa}\|x-a\|
\end{eqnarray*}
(thanks to $f(a)=0$ by \eqref{3.24}). Then \eqref{4.22} gives that
$$
(\beta-\frac{\epsilon}{\kappa})\|x-a\|\leq\tau f_+(x).
$$
By taking limits as $\beta\uparrow 1$, one gets
$$
\mathbf{d}(x, \mathcal{S})\leq\frac{\kappa\tau}{\kappa-\epsilon} f_+(x).
$$
This means that the convex-composite  inequality  system \eqref{3.21} has a local error bound with constant $\frac{\kappa\tau}{\kappa-\epsilon} >0$. The proof is complete.\hfill$\Box$

\medskip

The following corollary is a consequence result of \cref{th3.4}. It  gives a precise estimate of the local error bound modulus. 

\begin{corollary}
Let $f:=g\circ\psi$ be such that ${\rm bd}(\mathbf{S}_g)\subseteq g^{-1}(0)$ and $\bar x\in  \mathcal{S}$ such that $\triangledown \psi(\bar x)$ is surjective. Denote
\begin{equation}\label{3.26a}
  \tau^*(f,\bar x):=\{\tau>0: \exists\ \delta>0 \ s.t. \ \eqref{4.19}\ {\it holds\ for \ all \ } x\in \mathcal{S}\cap \mathbb{B}(\bar x, \delta) \}.
\end{equation}
Then 
\begin{equation}\label{3.27a}
  \tau(f,\bar x)=\tau^*(f,\bar x).
\end{equation}
\end{corollary}

{\bf Proof.} Consider the case that $\tau(f,\bar x)<+\infty$. Let $\tau\in (0, \tau(f,\bar x))$. Then there exists $\delta>0$ such that \eqref{4.8a} holds. Using the proof of the necessity part in \cref{th3.4}, one can verify that
$$
\tau^*(f,\bar x)\leq \tau
$$
and by letting $\tau\uparrow \tau(f,\bar x)$, one has
$$
\tau^*(f,\bar x)\leq \tau(f,\bar x).
$$

On the other hand, let $\tau\in (0, \tau^*(f,\bar x))$. Then for any $\epsilon>0$ sufficiently small,  by the proof of the sufficiency part in \cref{th3.4}, one has
$$
\tau(f,\bar x)\leq \frac{\kappa\tau}{\kappa-\epsilon}\rightarrow\tau\ ({\rm as} \ \epsilon\rightarrow 0^+)
$$
and thus 
$$
\tau(f,\bar x)\leq \tau.
$$
By letting $\tau\uparrow \tau^*(f,\bar x)$, one has
$$
\tau(f,\bar x)\leq \tau^*(f,\bar x).
$$
Hence \eqref{3.27a} holds. If $\tau(f,\bar x)=+\infty$, then we claim that $\tau^*(f,\bar x)=+\infty$ (otherwise, one can verify that $\tau(f,\bar x)<+\infty$ by the proof of the sufficiency part in \cref{th3.4}, a contradiction). The proof is complete.\hfill$\Box$

\medskip

\begin{theorem}\label{th3.5}
Let $f:=g\circ\psi$ be such that ${\rm bd}(\mathbf{S}_g)\subseteq g^{-1}(0)$ and $\bar x\in  \mathcal{S}$ such that $\triangledown \psi(\bar x)$ is surjective. Then convex-composite inequality \eqref{3.21} has a local error bound at $\bar x$ if and only if the convex inequality $g(y)\leq 0$ has a local error bound at $\psi(\bar x)$.
\end{theorem}

{\bf Proof.} Since $\triangledown \psi(\bar x)$ is surjective, it follows from \cite[Theorem 1.57]{Mordukhovich} that there exist $\kappa, \delta_0>0$ such that
\begin{equation}\label{4.28}
  \mathbf{d}(x, \psi^{-1}(y))\leq\kappa \|y-\psi(x)\|\ \ \forall (x,y)\in \mathbb{B}(\bar x, \delta_0)\times \mathbb{B}(\psi(\bar x), \delta_0).
\end{equation}
Note that $\mathcal{S}=\psi^{-1}(\mathbf{S}_g)$ and by applying \cite[Lemma 4.2]{WTY2023}, there are $\kappa, L>0$ and $r_0\in (0, \delta_0)$ such that \eqref{4.30} holds for all $x\in \mathbb{B}(\bar x,r_0)\cap \mathcal{S}$. 


{\it The necessity part.} By \cref{th3.4}, there exist $\tau,\delta\in (0,+\infty)$ such that \eqref{4.19} holds. We next show that there exist $\tau_1, \delta_1>0$ such that
\begin{equation}\label{4.26}
    N(\mathbf{S}_g, y)\cap \mathbb{B}_{\mathbb{Y}^*}\subseteq [0,\tau_1] \partial g(y),\ \ \forall y\in \mathbf{S}_g\cap \mathbb{B}(\bar y, \delta_1).
\end{equation}
Granting this, it follows from \cite[Theorem 2.2]{ZN2004} that  the convex  inequality $g(y)\leq 0$ has a local error bound at $\psi(\bar x)$.

Take $\delta_1\in (0,r_1)$ such that $\kappa\delta_1<\delta$. Let $y\in \mathbf{S}_g\cap \mathbb{B}(\psi(\bar x),\delta_1)$ and $y^*\in N(\mathbf{S}_g, y)\cap \mathbb{B}_{\mathbb{Y}^*}$. Then by \eqref{4.28}, one has
$$
\mathbf{d}(\bar x, \psi^{-1}(y))\leq\kappa\|\psi(\bar x)-y\|<\kappa\delta_1
$$
and so there is $x\in \mathbb{B}(\bar x, \kappa\delta_1)\subseteq \mathbb{B}(\bar x,\delta)$ such that $y=\psi(x)$. By virtue of \eqref{4.30}, one has
\begin{eqnarray*}
\triangledown \psi(x)^*(y^*)\in \triangledown \psi(x)^*(N(\mathbf{S}_g,\psi(x))\cap \mathbb{B}_{\mathbb{Y}^*})\subseteq\hat N(\mathcal{S},x)\cap L\mathbb{B}_{\mathbb{X}^*}
\end{eqnarray*}
and thus \eqref{4.19} implies that
\begin{eqnarray*}
\triangledown \psi(x)^*\Big(\frac{y^*}{L}\Big)\in\hat N(\mathcal{S},x)\cap \mathbb{B}_{\mathbb{X}^*}&\subseteq&  [0,\tau]\hat\partial f(x)\\
&=&  [0,\tau]\triangledown \psi(x)^*(\partial g(\psi(x)))\\
&=&\triangledown \psi(x)^*\big([0,\tau]\partial g(\psi(x))\big).
\end{eqnarray*}
Since $\triangledown \psi(x)^*$ is one-to-one, it follows that
$$
y^*\in [0, L\tau]\partial g(\psi(x)).
$$
This means that  \eqref{4.26} holds with $\tau_1:=L\tau$ and $\delta_1>0$.

{\it The sufficiency part}. Denote $\bar y:=\psi(\bar x)$. Then by virtue of \cite[Theorem 2.2]{ZN2004}, there exist $\tau_1,\delta_1>0$ such that \eqref{4.26} holds. Based on \cref{th3.4}, it suffices to show that there exist $\tau,\delta>0$ such that \eqref{4.19} holds.

Take $\delta>0$ such that $\psi(\mathbb{B}(\bar x,\delta))\subseteq \mathbb{B}(\psi(\bar x),\delta_1)$. Let $x\in \mathcal{S}\cap \mathbb{B}(\bar x,\delta)$ and $x^*\in\hat N(\mathcal{S}, x)\cap \mathbb{B}_{\mathbb{X}^*}$. By virtue of \eqref{4.30}, there is $y^*\in N(\mathbf{S}_g,\psi(x))\cap \mathbb{B}_{\mathbb{Y}^*}$ such that
\begin{equation}\label{4.31}
  \kappa x^*=\triangledown \psi(x)^*(y^*).
\end{equation}
By virtue of \eqref{4.26} and \cite[Lemma 4.1]{WTY2023}, one has 
\begin{eqnarray*}
\triangledown \psi(x)^*(y^*)&\in& \triangledown \psi(x)^*\Big( [0,\tau_1]\partial g(\psi(x))\Big)\\
 &=&  [0,\tau_1]\triangledown \psi(x)^*\Big(\partial g(\psi(x))\Big)\\
 &=& [0,\tau_1]\hat\partial f(x).
\end{eqnarray*}
This and \eqref{4.31} imply that 
$$
x^*\in [0,\frac{\tau_1}{\kappa}]\hat\partial f(x)
$$
and consequently \eqref{4.19} holds with $\tau:=\frac{\tau_1}{\kappa}$ and $\delta>0$. The proof is complete.\hfill$\Box$

\section{Conclusions}

This paper deals with error bounds of a general inequality defined by a proper lower semicontinuous function on an Asplund space. It is  described by  the subdifferentials of the given function at points inside the solution set.  Although  in absence of convexity assumption the dual characterization result on error bounds for the general inequality may not hold, 
 several necessary dual conditions for error bounds may  still  be obtained  in terms of Fr\'echet / Mordukhovich subdifferentials. Further, such conditions are also proved to be sufficient for error bounds of a  convex-composite inequality. Our work is an extension of  dual  characterizations  of error bounds via subdifferentials to the possible non-convex case.

\bibliography{WTY-2023}

\bibliographystyle{plain}

\end{document}